\documentclass[12pt, a4paper]{article}

\usepackage{mathptmx}
\usepackage{amsmath}
\usepackage{amsthm}
\usepackage[colorlinks,citecolor=blue,urlcolor=blue,linkcolor=blue,pagecolor=blue,linktocpage=true]{hyperref}
\usepackage{amssymb}
\usepackage{graphicx}
\usepackage[tight]{subfigure}
\usepackage{algorithmic}
\usepackage{algorithm} 

\usepackage{alphabeta}

\newcommand{\logn}{\ln}
\newcommand{\target}{\tau}
\newcommand{\error}{\varepsilon}
\newcommand{\lowerqu}{\lambda}
\newcommand{\upperqu}{\mu}
\newcommand{\lowerun}{\tilde\lambda}
\newcommand{\upperun}{\tilde\mu}

\newcommand{\cond}{\,|\,}\DeclareMathOperator {\E}{E}
\newcommand{\nbd}{B}
\newcommand{\erroralgo}{{\textsf{\textepsilon}}}
\newcommand{\lowerqualgo}{{\textsf{\textlambda}}}
\newcommand{\lowerunalgo}{{\tilde{\textsf{\textlambda}}}}
\newcommand{\Nalgo}{{\textsf{N}}}
\newcommand{\salgo}{{\textsf{s}}}
\newcommand{\kalgo}{{\textsf{k}}}

\newtheorem{proposition}{Proposition}
\newtheorem{theorem}{Theorem}
\theoremstyle{definition} 

\begin{document}

\title{
Simulating a coin with irrational bias\\ using rational arithmetic
}

\author{Luis Mendo
\\
Universidad Polit\'ecnica de Madrid\\
\texttt{luis.mendo@upm.es}
}


\maketitle

\begin{abstract}
An algorithm is presented that, taking independent Bernoulli random variables with parameter $1/2$ as inputs and using only rational arithmetic, simulates a Bernoulli random variable with possibly irrational parameter $\target$. It requires a series representation of $\target$ with positive, rational terms, and a rational bound on its truncation error that converges to $0$. The number of required inputs has an exponentially bounded tail, and its mean is at most $3$. The number of arithmetic operations has a tail that can be bounded in terms of the sequence of truncation error bounds.

The algorithm is applied to two specific values of $\target$, including Euler's constant, for which obtaining a simple simulation algorithm was an open problem.


\end{abstract}

\section{Introduction}
\label{part: intro}

Consider the problem of generating a Bernoulli random variable $Y$ with known parameter $\target \in (0, 1)$. The only source of randomness is a sequence of independent, identically distributed Bernoulli variables $X_k$ with parameter $1/2$, and only rational-number computations can be performed. A \emph{sequential} algorithm will be used, whereby the number of consumed inputs is random, and is governed by a certain \emph{stopping rule}.

The enounced problem has applications in random number generation and simulation. In the following, independent Bernoulli random variables with parameter $p$ will be referred to as ``$p$-coins'', or ``unbiased coins'' if $p=1/2$. The restriction to use unbiased coins as inputs is a natural one, as these constitute the most ``basic'' form of randomness. Unbiased coins can easily be obtained from $p$-coins, even if $p$ is unknown, using the well known von Neumann's procedure \cite{VonNeumann51} or its refinement by Peres \cite{Peres92}. Requiring only rational arithmetic operations has obvious advantages in terms of simplicity and precision. Clearly, rational arithmetic can be reduced to integer arithmetic by storing each rational number $x = n/d$ as $n$ and $d$ separately, and implementing operations on $x = n/d$ and $x' = n'/d'$ by means of integer arithmetic with $n$, $d'$, $n$ and $d'$.

The stated problem is easy when $\target$ is a rational number $n/d$. Let $b$ be the number of digits in the binary expression of $d$. Then it suffices to sample $b$ unbiased coins, interpret the result as an integer $t \in \{0,1,\ldots,2^b-1\}$, repeat until $t \leq d-1$, and then output $1$ if $t \leq n-1$ or $0$ otherwise. The general problem when $\target$ is not necessarily rational is more interesting.

An algorithm for producing a $\target$-coin from unbiased coins is termed a \emph{Buffon machine} in \cite{Flajolet11}.
This is related to the \emph{Bernoulli factory} problem \cite{Keane94}, which consists in generating an $f(p)$-coin from $p$-coins, when the function $f$ is known and $p$ is unknown. This problem has been extensively studied, and it is known how characteristics of algorithms that can solve it are restricted by properties of the function \cite{Keane94, Nacu05}. Efficient algorithms are known for several classes of functions \cite{Latuszynski11, Huber16, Mendo19}.

This paper describes a general algorithm that solves the above problem when there exists a representation of $\target$ as a series with rational terms. The algorithm is described in \S \ref{part: descr, basic prop}, and its basic properties are addressed. The complexity of the algorithm is analysed in \S \ref{part: compl}. Application to specific values, including Euler's constant $\gamma$ and $\pi/4$, is discussed in \S \ref{part: applic}. Conclusions and future work are presented in \S \ref{part: concl}. Proofs to all results are given in \S \ref{part: proofs}.

The following notation and definitions are used. A random variable $V$ is referred to as a shifted geometric variable if $V-1$ is geometric. Given two non-negative functions $f$ and $g$, $f(x)$ is said to be $O(g(x))$ if there exist $K$ and $x_0$ such that $f(x) \leq K g(x)$ for all $x \geq x_0$. The natural and binary logarithms of $x$ are written as $\logn x$ and $\log_2 x$ respectively. The number of digits in the binary expression of a positive integer $t$ is denoted as $\nbd(t)$. For $n \in \mathbb N$, the notation $n!!$ represents the double factorial, that is, $n (n-2) \cdots 3 \cdot 1$ for $n$ odd and $n (n-2) \cdots 4 \cdot 2$ for $n$ even.

\section{Algorithm description and basic properties}
\label{part: descr, basic prop}

The proposed algorithm is inspired by \cite[algorithm 2]{Latuszynski11}, which uses a continuous uniform random variable $U$ on $(0,1)$, a sequence of monotonically increasing lower bounds $\lowerqu_k$ that converge to $\target$, and a sequence of monotonically decreasing upper bounds $\upperqu_k$ that converge to $\target$. The algorithm in the cited reference generates $U$ and then computes $\lowerqu_k$, $\upperqu_k$ for successive $k$ until $U \leq \lowerqu_k$ or $U > \upperqu_k$ (this occurs for finite $k$ with probability $1$). The output $Y$ is then $1$ if $U \leq \lowerqu_k$ and $0$ if $U > \upperqu_k$ for the last $k$. 

Generating a \emph{continuous} uniform random variable in a computer simulation poses numerical precision problems, and is not compatible with using \emph{rational} arithmetic.
The reason is that computer simulations typically represent non-integer numerical values using floating-point data types, which implies that the number of significant digits that can be used is limited to a fixed value.
Instead, the algorithm to be presented does not generate $U$ explicitly, but only the information about it which is \emph{necessary} at each iteration $k$ to determine if $U$ is below $\lowerqu_k$, above $\upperqu_k$, or between the two bounds.
This avoids the loss of accuracy that would be incurred when trying to represent the exact value of $U$.

Let $\target \in (0,1)$ be represented as a convergent series with positive, rational terms $a_j$:
\begin{equation}
\label{eq: serie, pos}
\target = \sum_{j=1}^\infty a_j.
\end{equation}
Series expressions of this form are available for the majority of commonly used constants. In addition, a bound $\error(N)$ for the truncation error is assumed to be known, and to be computable with operations involving only rational numbers:
\begin{equation}
\label{eq: serie, error}
\target - \sum_{j=1}^N a_j \leq \error(N),
\end{equation}
where $\error(N)$ is monotonically non-increasing with $\lim_{N \rightarrow \infty} \error(N) = 0$. The function $\error$ will be referred to as \emph{error function}.

An alternating series
\begin{equation}
\label{eq: serie, alt}
\target = \sum_{j=1}^\infty (-1)^{j+1} b_j
\end{equation}
with terms that decrease monotonically in absolute value can be rewritten in the form \eqref{eq: serie, pos} with $a_j = b_{2j-1}-b_{2j}$. If $\lim_{j \rightarrow \infty} b_j = 0$, Leibniz's rule \cite[theorem 10.14]{Apostol67} implies that the series converges, and $\target - \sum_{j=1}^N a_j \leq b_{2N+1}$. Thus a simple characterization of the truncation error is possible in this case, namely $\error(N) = b_{2N+1}$. Series with negative terms or with alternating signs opposite from those in \eqref{eq: serie, alt} are reduced to the preceding cases by considering $1-\target$ instead of $\target$; and then it suffices to replace the algorithm output $Y$ by $Y'=1-Y$ to achieve $\Pr[Y'=1]=\target$.

The monotonicity requirement for the error function does not impose any restriction, because any error function can be modified to fulfil this condition, simply replacing $\error(n)$ by its cumulative minimum $\error'(n) = \min\{\error(1),\ldots,\error(n)\}$. This can be done because, since the series has positive terms, the error bound $\error(i)$ is valid not only for $\sum_{j=1}^i a_j$, but also for any $\sum_{j=1}^n a_j$ with $n > i$. As will be seen, the algorithm to be presented uses the error bound for the sum with $n$ terms after it has already used the error bound for the sum with $n-1$ terms. Thus the cumulative minimum $\error'(n)$ can be efficiently obtained as
\begin{equation}
\error'(n) = \begin{cases}
\error(n) & \text{if } n  =1 \\
\min\{\error'(n-1), \error(n)\} & \text{if } n \geq 2.
\end{cases}
\end{equation}

The proposed algorithm consists of a random number $M$ of \emph{iterations}. At the beginning of iteration $k$, the continuous uniform variable $U$ is known to be in an interval $(\lowerqu_{k-1}, \upperqu_{k-1}]$ resulting from the previous iteration. The iteration \emph{shrinks} this interval to a new interval $(\lowerqu_k, \upperqu_k] \subset (\lowerqu_k-1, \upperqu_k-1]$. These intervals are \emph{quantized}, with finer resolution as the algorithm progresses. Specifically, the endpoints of the interval $(\lowerqu_k, \upperqu_k]$ are multiples of $2^{-(k+1)}$, and the length of the interval is $2^{-k}$. Thus each quantized interval is half as wide as the one from the preceding iteration. The shrinking and quantizing conditions leave only three possible choices for the interval $(\lowerqu_k, \upperqu_k]$ given $(\lowerqu_{k-1}, \upperqu_{k-1}]$:
\begin{align}
\label{eq: choice of s k}
\lowerqu_k = \lowerqu_{k-1} + s_k\cdot 2^{-(k+1)},\quad \upperqu_k = \lowerqu_k + 2^{-k},\quad s_k \in \{0, 1, 2\}.
\end{align}
Thus $(\lowerqu_k, \upperqu_k]$ is the lower half, the middle half or the upper half of $(\lowerqu_{k-1}, \upperqu_{k-1}]$ for $s_k = 0$, $1$ or $2$ respectively (see Figure \ref{fig: from iteration k-1 to k} below).

The choice of $s_k$ is dictated by intermediate, unquantized bounds $\lowerun_k$ and $\upperun_k$ computed from the series representation of $\target$:
\begin{align}
\label{eq: lower bound un k}
\lowerun_k &= \sum_{j=1}^{N_k} a_j, \\
\label{eq: upper bound un k}
\upperun_k &= \lowerun_k + \error(N_k),
\end{align}
which define an unquantized interval $(\lowerun_k, \upperun_k]$ that also shrinks at each iteration. More specifically, knowing $N_{k-1}$, $\lowerun_{k-1}$ and $\upperun_{k-1}$ from the previous iteration, the new bounds $\lowerun_k$ and $\upperun_k$ are obtained adding series terms up to a certain index $N_k \geq N_{k-1}$:
\begin{equation}
\label{eq: lower bound un k sep}
\lowerun_k = \lowerun_{k-1} + \sum_{j=N_{k-1}+1}^{N_k} a_j, \\
\end{equation}
and computing the corresponding truncation error $\error(N_k)$ to be used in \eqref{eq: upper bound un k}. The number of terms $N_k$ at iteration $k$ is chosen as the smallest value such that $(\lowerun_k, \upperun_k] \cap (\lowerqu_{k-1}, \upperqu_{k-1}]$ is contained in one of the three possible quantized intervals $(\lowerqu_k, \upperqu_k]$ defined by \eqref{eq: choice of s k}, which also determines the choice of $s_k$. The reason for this is that if $\target$ is in $(\lowerun_k, \upperun_k]$ and this interval is contained in one of the three quantized intervals, $\target$ is assured to be contained in that quantized interval. Besides, if part of the interval $(\lowerun_k, \upperun_k]$ exceeds the boundaries of $(\lowerqu_{k-1}, \upperqu_{k-1}]$ that part can be disregarded (only the intersection matters), because $\target$ is known not to be outside $(\lowerqu_{k-1}, \upperqu_{k-1}]$. Note that both sequences $N_k$ and $s_k$ are deterministic. Figure \ref{fig: from iteration k-1 to k} illustrates the steps involved in moving from $(\lowerqu_{k-1}, \upperqu_{k-1}]$ to $(\lowerqu_k, \upperqu_k]$.

\begin{figure}%
\centering%
\subfigure[Case $(\lowerun_k, \upperun_k{]} \subset (\lowerqu_{k-1}, \upperqu_{k-1}{]}$]{%
\label{fig: fig: from iteration k-1 to k, subset}%
\includegraphics[width=.9\textwidth]{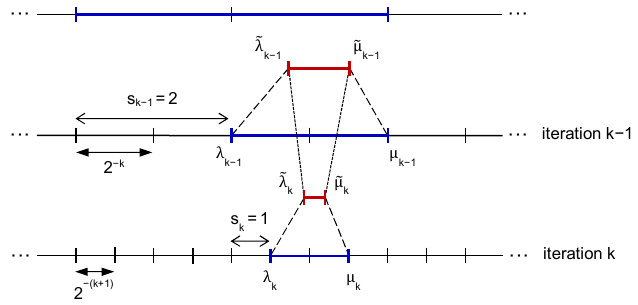}%
}\\[5mm]
\subfigure[Case $(\lowerun_k, \upperun_k{]} \not\subset (\lowerqu_{k-1}, \upperqu_{k-1}{]}$]{%
\label{fig: from iteration k-1 to k, no subset}%
\includegraphics[width=.9\textwidth]{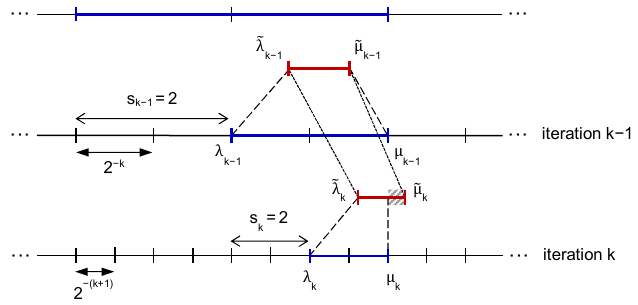}%
}%
\caption{Steps in moving from iteration $k-1$ to $k$
}%
\label{fig: from iteration k-1 to k}%
\end{figure}%

According to the above, the rules for selecting $N_k$ and $s_k$ at iteration $k$ are:
\begin{equation}
\label{eq: rules for N k and s k}
\begin{split}
\text{if }\, &\upperun_k \leq \lowerqu_{k-1} + 2^{-k} \Rightarrow s_k=0;\\
\text{elseif }\, &\lowerun_k > \lowerqu_{k-1} + 2^{-k} \Rightarrow s_k=2;\\
\text{elseif }\, &\lowerun_k > \lowerqu_{k-1} + \frac 1 2 \cdot 2^{-k} \text{ and } \upperun_k \leq \lowerqu_{k-1} + \frac 3 2 \cdot 2^{-k} \Rightarrow s_k=1;\\
\text{else }\, & \text{a larger } N_k \text{ is needed}.
\end{split}
\end{equation}
Thus, starting from the partial sum $\lowerun_{k-1}$ from the previous iteration, which contains $N_{k-1}$ terms, new terms are added one by one, computing tentative values of $\lowerun_k$ and $\upperun_k$ for each new partial sum, checking conditions \eqref{eq: rules for N k and s k}, and stopping as soon as one of the three conditions holds.

Once the new quantized interval $(\lowerqu_k, \upperqu_k]$ has been obtained in iteration $k$, a random input $X_k$ is used to decide if $U$ is in $(\lowerqu_k, \upperqu_k]$. Specifically, the algorithm uses initial values $\lowerun_0 = 0$, $\upperun_0 = 1$ and $(\lowerqu_0, \upperqu_0] = (0, 1]$. At this point, $U$ is only known to be in this interval.
The algorithm proceeds with the first iteration, $k=1$, and computes the interval $(\lowerqu_1, \upperqu_1]$, which has length $1/2$. Therefore $\Pr[U \in (\lowerqu_1, \upperqu_1] \cond U \in (0, 1]] = \Pr[U \in (\lowerqu_1, \upperqu_1]] = 1/2$. Thus an input $X_1$ is taken to randomly \emph{decide} if $U$ is in $(\lowerqu_1, \upperqu_1]$. If it is, the algorithm moves to iteration $k=2$. In this iteration $U$ is known to be in $(\lowerqu_1, \upperqu_1]$, and the distribution of $U$ \emph{conditioned} on this information is uniform on that interval. The new interval $(\lowerqu_2, \upperqu_2]$ is a subset of $(\lowerqu_1, \upperqu_1]$ with half its length, thus again $\Pr[U \in (\lowerqu_2, \upperqu_2] \cond U \in (\lowerqu_1, \upperqu_1]] = 1/2$, which can be simulated with a new random input $X_2$. This way, knowledge about $U$ is refined at each iteration, reducing the interval in which $U$ is known to be to a new interval with half the length, which also contains $\target$.

Eventually (with probability $1$) there will be one iteration, with index $M$, for which the random input $X_M$ indicates that $U \notin (\lowerqu_M, \upperqu_M]$. This is the last iteration. At this point $s_M$ is known, and it is also known that $\target \in (\lowerqu_M, \upperqu_M]$. If $s_M = 0$ the interval $(\lowerqu_M, \upperqu_M]$ is the lower half of  $(\lowerqu_{M-1}, \upperqu_{M-1}]$, thus $U > \upperqu_M \geq \target$, and the output $Y$ is $0$. Similarly, if $s_M=2$ the output is $Y=1$, because $U \leq \lowerqu_M < \target$. If $s_M=1$ the interval $(\lowerqu_M, \upperqu_M]$ is in the middle of $(\lowerqu_{M-1}, \upperqu_{M-1}]$ and a final input $X_{M+1}$ is needed to decide if $U$ is in the lower or upper quarter of $(\lowerqu_{M-1}, \upperqu_{M-1}]$, both events being equally likely, to determine if the output $Y$ is $1$ or $0$ respectively.

The number of iterations $M$ is, by construction, a shifted geometric random variable with parameter $1/2$, and thus for any $m \in \mathbb N$
\begin{equation}
\label{eq: Pr M geq}
\Pr[M \geq m] = 2^{-m+1}.
\end{equation}
Denoting the total number of required inputs by $L$, it is clear that $L=M+1$ or $L=M$, depending on whether one last input is needed to generate the output or not. 

As is apparent from the foregoing description, although the continuous variable $U$ is helpful for explaining the process, its exact value is not actually needed, and is never generated. $U$ is only known to be in intervals of decreasing size, that are constructed by the algorithm based on the partial sums and error bounds of the series.

The described procedure can be compared with that given in \cite[section 1]{Flajolet11} to simulate a rational constant $\target \in (0, 1)$: using the binary representation of $\target$, which in the rational case is completely known (it is either finite or repeating), output the $i$-th binary digit where $i$ is given by a shifted geometric random variable with parameter $1/2$. To extend this approach for irrational $\target$ its binary representation, which is not fully known, would have to be computed up to the $i$-th digit. That is similar to what the algorithm presented here does: it computes increasingly accurate approximations of $\target$ until a sufficient level of accuracy, given by the shifted geometric random variable $M$, is achieved.

Another related approach is the ``interval algorithm'' in \cite{Han97}, which in general transforms a discrete, finite-support distribution into another, both distributions being known. Particularized to an unbiased coin as input and a biased coin with parameter $\target$ as output, the referred algorithm iteratively replaces an interval, initialised as $[0,1)$, by either its lower or upper half, as indicated by a random input, until the interval no longer contains $\target$. The main differences with the algorithm described in this paper are that \cite{Han97} uses disjoint subintervals and assumes $\target$ to be known exactly. Here, in contrast, overlapping subintervals are used, and $\target$ needs not be known exactly (which allows using rational arithmetic).

Based on the above the algorithm can be precisely stated; see Algorithm \ref{algo: main}. Its inputs are:
\begin{enumerate}
\item
A series with positive terms $a_j$, as given by \eqref{eq: serie, pos}, that converges to the target value $\target$.
\item
A function that computes a bound $\error(N)$ for the error of approximating the series with the first $N$ terms, as given by \eqref{eq: serie, error}, with $\lim_{N \rightarrow \infty} \error(N) = 0$. If necessary, its cumulative minimum should be taken in order to make the function monotonically non-increasing.
\item
A sequence of independent Bernoulli random variables $X_k$ with parameter $1/2$.
\end{enumerate}
The output is a Bernoulli random variable $Y$ with parameter $\target$.

\begin{algorithm}
\caption{Simulation of a constant using unbiased coins}
\label{algo: main}
\begin{algorithmic}
\STATE $\Nalgo \leftarrow 0$, $\lowerunalgo \leftarrow 0$, $\erroralgo \leftarrow 1$
\STATE $\lowerqualgo \leftarrow 0$, $\salgo \leftarrow 0$, $\kalgo \leftarrow 0$
\REPEAT
	\STATE $\kalgo \leftarrow \kalgo+1$
	\STATE $\lowerqualgo \leftarrow \lowerqualgo + \salgo\cdot 2^{-\kalgo}$
	\WHILE{\NOT ($\lowerunalgo + \erroralgo \leq \lowerqualgo + 2^{-\kalgo}$ \OR $\lowerunalgo > \lowerqualgo + 2^{-\kalgo}$ \OR \\ ($\lowerunalgo > \lowerqualgo + (1/2) \cdot 2^{-\kalgo}$ \AND  $\lowerunalgo + \erroralgo \leq \lowerqualgo + (3/2) \cdot 2^{-\kalgo}$))}
		\STATE $\Nalgo \leftarrow \Nalgo+1$
		\STATE $\lowerunalgo \leftarrow \lowerunalgo + a_\Nalgo$
		\STATE $\erroralgo \leftarrow \error(\Nalgo)$
	\ENDWHILE
	\IF{$\lowerunalgo + \erroralgo \leq \lowerqualgo + 2^{-\kalgo}$}
		\STATE $\salgo \leftarrow 0$
	\ELSIF{$\lowerunalgo > \lowerqualgo + 2^{-\kalgo}$}
		\STATE $\salgo \leftarrow 2$
	\ELSE 
		\STATE $\salgo \leftarrow 1$
	\ENDIF
\UNTIL $X_\kalgo = 0$
\IF{$\salgo = 0$}
	\STATE $Y \leftarrow 0$
\ELSIF{$\salgo = 2$}
	\STATE $Y \leftarrow 1$
\ELSE 
	\STATE $Y \leftarrow X_{\kalgo+1}$
	\ENDIF
\end{algorithmic}
\end{algorithm}

Algorithm \ref{algo: main} is seen to consist of three parts: initialisation of variables (first two lines), iterations (main loop: ``repeat''), and output generation (final ``if'' block). The inner loop (``while'') updates the partial sum of the series and its error bound. A difference from the description in the preceding paragraphs is that sequences of variables such as $N_k$, $\lowerun_k$, $\lowerqu_k$, $s_k$, which were indexed by the iteration number $k$, are here expressed more compactly by single variables $\Nalgo$, $\lowerunalgo$, $\lowerqualgo$, and $\salgo$ that are \emph{updated} at each iteration. Similarly, the variable $\erroralgo$ stores the value of the latest error bound, $\error(N_k)$; and $\kalgo$ is the current iteration index (which defines resolution).

\begin{theorem}[Basic properties]
\label{theo: basic}
Algorithm \ref{algo: main} satisfies the following:
\begin{enumerate}
\item
The algorithm terminates with probability $1$.
\item
The output $Y$ is a Bernoulli random variable with parameter $\target$.
\item
If each term $a_j$ in the series \eqref{eq: serie, pos} and the error bound $\error(N)$ defined by \eqref{eq: serie, error} can be computed with a finite number of operations involving rational numbers, the algorithm can be implemented using rational arithmetic.
\end{enumerate}
\end{theorem}

\section{Complexity analysis}
\label{part: compl}

The complexity of a sequential algorithm is primarily determined by the \emph{number of required inputs}, because consuming a new input is typically considered more costly than the arithmetical operations needed to process it. However, it is also important to analyse the \emph{number of required arithmetical operations}, to ensure that it is not unrealistically large. In this regard, assessing its order of magnitude is usually enough.
 
Algorithm \ref{algo: main} consumes one input for each iteration, and possibly one additional input to generate the output. As for arithmetical operations, the algorithm uses two groups thereof, from the point of view of a complexity analysis:
\begin{enumerate}
\item
Operations that are carried out to obtain each new term of the series, and to update the partial sum and error bound. The total number of these operations is roughly proportional to the total number of terms in the partial sum when the algorithm ends, $N_M$.
\item
Operations that are needed to update the variables used by the algorithm. The total number of these is proportional to the number of iterations, $M$.
\end{enumerate}

From the preceding it is clear that a good characterization of algorithm complexity can be obtained by analysing the number of iterations $M$, the number of inputs $L$, and the number of series terms when the algorithm terminates, $N_M$. The first two are characterized very easily, because the distribution of $M$ is known.

An algorithm that uses $L$ inputs is \emph{fast}, as defined by Nacu and Peres \cite{Nacu05}, if $L$ is exponentially bounded, that is, if there exist $C>0$, $\rho < 1$ such that for all $l \in \mathbb N$
\begin{equation}
\label{eq: fast}
\Pr[L > l] \leq C \rho^l.
\end{equation}

\begin{theorem}[Number of inputs]
\label{theo: number of inputs}
In Algorithm \ref{algo: main}, the number of required inputs satisfies
\begin{equation}
\Pr[L>l] \leq 2^{-l+1}
\end{equation}
for $l \in \mathbb N$, and thus the algorithm is fast in the sense of Nacu-Peres. In addition,
\begin{equation}
\label{eq: bound on E L}
2 \leq \E[L] \leq 3.
\end{equation}
\end{theorem}

Thus, according to this theorem, the number of inputs used by the algorithm has an exponentially bounded tail and its average value is very small. Indeed, \cite[theorem 6]{Kozen14} establishes that any algorithm that outputs a Bernoulli random variable with parameter $\target$ from inputs with parameter $1/2$ must use at least $2$ inputs on average, except when $\target$ is a dyadic number. Therefore the average number of inputs consumed by Algorithm \ref{algo: main} is close to the optimum. Moreover, even if the number of required inputs in a given realization of the algorithm can potentially be much larger than its average value, the probability that this happens is very small thanks to the exponential-bound property.

The distribution of $N_M$ is more difficult to characterize, because it depends on the error function $\error$: the more slowly this function decreases, the more likely it is for $N_M$ to take larger values. However, as has been discussed, $N_M$ only affects the number of required arithmetical operations and thus it suffices to know its order or magnitude. The following result establishes a bound on $\Pr[N_M > n]$, and gives a sufficient condition for $\E[N_M]$ to be finite.

\begin{theorem}[Number of series terms]
\label{theo: number of series terms}
The number of series terms used by Algorithm \ref{algo: main} satisfies
\begin{equation}
\Pr[N_M > n] < 4 \error(n)
\end{equation}
for $n \in \mathbb N$, where $\error$ is the error function. In addition, if $\error(n)$ is $O(1/n^r)$ for some $r>1$, $\E[N_M]$ is finite. 
\end{theorem}

This theorem describes how the truncation error bound influences the number of required arithmetical operations; namely, $\Pr[N_M > n]$ is $O(\error(n))$. Furthermore, $\error(n)$ asymptotically decreasing as the inverse of a power with exponent greater than $1$ is sufficient to ensure that $\E[N_M]$ is finite. This requirement on $\error$ is not very stringent. The two detailed examples to be presented in \S \ref{part: applic} will satisfy this condition (and $\E[N_M]$ will be seen to be not only finite but very small).

A conceivable modification of Algorithm \ref{algo: main} would be to only allow cases $s_k = 0$ and $s_k = 2$ in each iteration (like \cite{Han97} does). With this approach the third condition in \eqref{eq: rules for N k and s k} is eliminated, and series terms are added until the interval $(\lowerun_k, \upperun_k] \cap (\lowerqu_{k-1}, \upperqu_{k-1}]$ is contained either in the lower half or in the upper half of $(\lowerqu_{k-1}, \upperqu_{k-1}]$. This way the final iteration never requires an additional input to generate the output. The problem with this method is that, depending on the value of $\target$, shrinking $(\lowerun_k, \upperun_k]$ until it fits into one of the two halves of the previous quantized interval may require an arbitrarily large number of series terms. Furthermore, if $\target$ is a dyadic number there is a non-zero probability that the algorithm does not terminate. Algorithm $1$ increases the average number of inputs from the optimum $2$ to at most $3$, but in return it terminates with probability $1$ and $\Pr[N_M > n]$ is bounded. (Note that the procedure in \cite{Han97} always terminates, but assumes perfect knowledge of $\target$).

\section{Application}
\label{part: applic}

Two specific cases will be considered in detail: Euler's constant $\gamma$ (\S \ref{part: sim gamma}), and $\pi/4$ (\S \ref{part: sim pidiv4}). These are in themselves interesting; especially the former, as simulating $\gamma$ without real-number arithmetic is one of the open problems mentioned in \cite{Flajolet11}. In addition, they illustrate the algorithm's performance in two different situations: a series of positive terms with error decaying as a power law, and an alternating series with exponential error decay. A few additional examples are then briefly discussed (\S \ref{part: other examples}).

\subsection{Simulation of $\gamma$}
\label{part: sim gamma}

Euler's constant is defined as $\gamma = \lim_{n \rightarrow \infty} (-\logn n + \sum_{i=1}^n 1/i) = 0.5772156\ldots$ Many series are known that converge to $\gamma$; see for example \cite{Addison67, Sondow10b, Blagouchine16a}.
The following one \cite{Addison67, Sondow10b} is of interest for application of Algorithm \ref{algo: main}:
\begin{equation}
\label{eq: series, gamma, orig}
\gamma = \frac 1 2 + \sum_{j=1}^\infty \frac{\nbd(j)}{2j(2j+1)(2j+2)} 
\end{equation}
where $\nbd(n)$ is the number of binary digits of the positive integer $n$, as defined in \S \ref{part: intro}. This can be computed using only integer operations, namely successive values $b=1,2,\ldots$ are tried, and the output is the first $b$ such that $2^b>n$.

The series \eqref{eq: series, gamma, orig} can be rewritten in the form \eqref{eq: serie, pos} with
\begin{equation}
\label{eq: series, gamma, aj}
a_j = \begin{cases}
1/2 & \text{if } j=1 \\
\displaystyle
\frac{\nbd(j-1)}{2j(2j-1)(2j-2)} & \text{if } j \geq 2.
\end{cases}
\end{equation}
A rational bound can easily be obtained for the truncation error of this series.

\begin{proposition}
\label{prop: error bound for gamma}
For $N \geq 2$, the series defined by \eqref{eq: series, gamma, aj} satisfies $\gamma - \sum_{j=1}^N a_j < \error(N)$ with
\begin{equation}
\label{eq: series, gamma, error, no monot}
\error(N) = \frac{2 + \nbd(N-1) + 1/(N-1)}{16 (N-1)^2}.
\end{equation}
\end{proposition}

This error bound clearly converges to $0$, but it is not monotonic, due to the term $\nbd(N-1)$ in the numerator. Namely, $\nbd(2^t) = \nbd(2^t-1)+1$ for any positive integer $t$, which can cause \eqref{eq: series, gamma, error, no monot} to increase when $N$ changes from $2^t$ to $2^t+1$. Indeed, it can be seen that $\error(N+1) > \error(N)$ for $N = 2^4, 2^5, 2^6, \ldots$ and $\error(N+1) < \error(N)$ otherwise. Nonetheless, as discussed in \S \ref{part: descr, basic prop}, monotonicity can be achieved by redefining the error function as
\begin{equation}
\label{eq: series, gamma, error}
\error(N) =
\begin{cases}
1/2 & \text{if } N=1 \\
\displaystyle
\min\left\{\error(N-1), \frac{2 + \nbd(N-1) + 1/(N-1)}{16 (N-1)^2} \right\} & \text{if } N \geq 2.
\end{cases}
\end{equation}
In addition, this function is easily shown to be $O(1/n^r)$ for any $r<2$.

Thus the series defined by \eqref{eq: series, gamma, aj} and its truncation error bound \eqref{eq: series, gamma, error} satisfy all the requirements in \S \ref{part: descr, basic prop} and \S \ref{part: compl} (series with positive rational terms that can be computed easily; error bound that is monotonically non-increasing, tends to $0$, and is $O(1/(n^r)$ for some $r>1$). Therefore Algorithm \ref{algo: main} can be applied to \eqref{eq: series, gamma, aj} and \eqref{eq: series, gamma, error}, and the results in Theorems \ref{theo: basic}--\ref{theo: number of series terms}
hold.

Table \ref{tab: gamma} shows the sample mean of $Y$, $L$ and $N_M$ obtained from running the algorithm $10^8$ times. The sample mean of $Y$ differs from $\gamma$ by $0.000027$. This difference is comparable to the standard deviation of the average of $10^8$ Bernoulli variables with parameter $\gamma$, which is $(\gamma(1-\gamma)/10^8)^{1/2} = 0.000049$. The average number of inputs is only slightly greater than $2$, and therefore very close to the optimum.

The average number of series terms is also seen to be very small. In view of \eqref{eq: series, gamma, aj} and \eqref{eq: series, gamma, error}, updating the partial sum of the series with a new term and computing the corresponding error bound requires around $20$ arithmetical operations (the exact number of operations depends on implementation details such as whether intermediate results are stored; note also that the first terms, which are needed with higher probability, require fewer operations). Updating the algorithm variables in each iteration adds a small amount of computational burden. Therefore a simulation requires an average number of arithmetical operations of the order of several tens.

\begin{table}
\caption{Sample results for simulation of $\gamma$}
\label{tab: gamma}
\begin{center}
  \begin{tabular}{|r|r|r|r|}
    \hline
    \multicolumn{1}{|p{2.7cm}}{\centering $\target$} & \multicolumn{1}{|p{2.7cm}}{\centering Average of $Y$} &  \multicolumn{1}{|p{2.7cm}}{\centering Average of $L$} & \multicolumn{1}{|p{2.7cm}|}{\centering Average of $N_M$} \\ \hline\hline
    $0.577215\ldots$ & $0.577243$ & $2.0250$ & $3.0053$ \\ \hline
  \end{tabular}
\end{center}
\end{table}

\begin{figure}
\centering
\includegraphics[width = .75\textwidth]{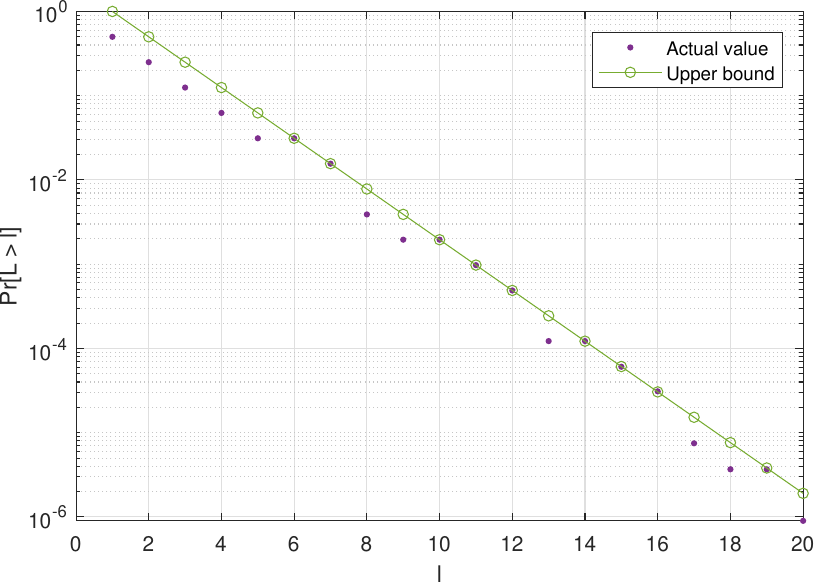}
\caption{
Sample estimate of $\Pr[L > l]$ for simulation of $\gamma$}
\label{fig: gamma_PrL}
\end{figure}

\begin{figure}
\centering
\includegraphics[width = .75\textwidth]{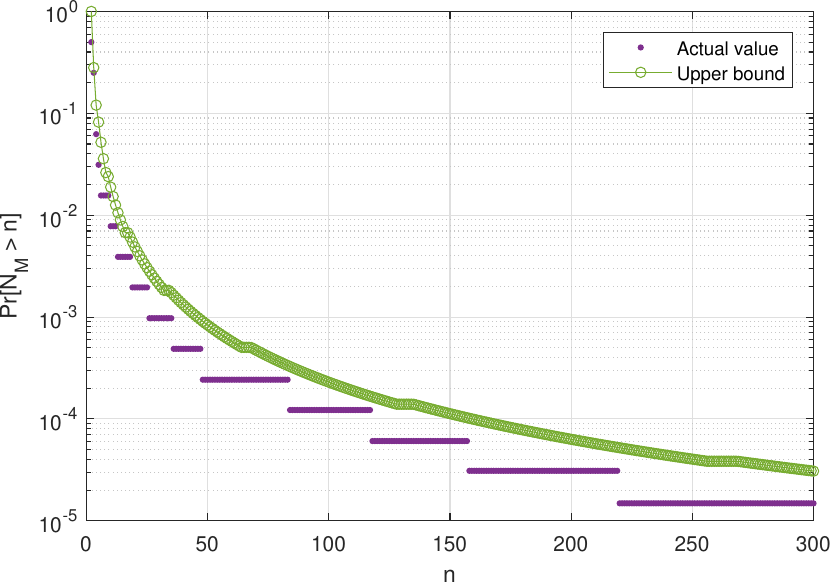}
\caption{
Sample estimate of $\Pr[N_M > n]$ for simulation of $\gamma$}
\label{fig: gamma_PrNM}
\end{figure}

Figures \ref{fig: gamma_PrL} and \ref{fig: gamma_PrNM} depict the sample estimates of $\Pr[L > l]$ and $\Pr[N_M > n]$, together with their respective bounds $2^{-l+1}$ and $4\error(n)$, for $l, n \in \mathbb N$. The actual value (filled circle) is always below or at the same height as the corresponding bound (empty circle). $\Pr[L > l]$ equals its bound whenever $s_l=1$, or half its bound otherwise. Observe how $\Pr[N_M > n]$ consists of runs of equal values, caused by the fact that $N_M$ cannot take any value between $N_k$ and $N_{k+1}$. Also, the bound of $\Pr[N_M > n]$ has short runs of equal values near all powers of $2$ except for the smallest ones, corresponding to an increase in the original bound \eqref{eq: series, gamma, error, no monot}.

Some optimization is possible if the algorithm is to be repeatedly applied to generate several independent values of $Y$. Namely, partial results can be stored to avoid computing them more than once. This applies to values of algorithm variables $\Nalgo, \lowerunalgo, \erroralgo, \lowerqualgo, \salgo$ at the end of each iteration. Thus a simulation only has to compute new values if it advances to an iteration index that has never been reached before.

This also applies to the function $\nbd$, even within a given simulation: once $\nbd(n)$ is known for some $n$, computing $\nbd(n')$ for $n'>n$ can start from the previous value to save operations. Thus in \eqref{eq: series, gamma, aj} and \eqref{eq: series, gamma, error} the total number of operations required by all evaluations of function $\nbd$ is that corresponding to its largest input argument.

The application of Algorithm \ref{algo: main} to $\gamma$ answers an open question posed in \cite[section 6]{Flajolet11}, namely devising a ``natural'' experiment with probability of success $\gamma$. The procedure presented here is ``natural'' in the sense that it only requires rational arithmetic and its complexity is very small, both in terms of consumed inputs and of number of operations.

\subsection{Simulation of $\pi/4$}
\label{part: sim pidiv4}

Rational multiples of $\pi$ are a classical example for the simulation of Bernoulli random variables. Consider Euler's Machin-like formula \cite{Nishiyama13},
\begin{equation}
\label{eq: Euler Machin}
\frac{\pi}{4} = \arctan \frac 1 2 + \arctan \frac 1 3.
\end{equation}
The Taylor expansion of the $\arctan$ function \cite[section 4.4]{Abramowitz72},
\begin{equation}
\arctan x = \sum_{j=1}^\infty \frac{(-1)^{j+1}}{2j-1}x^{2j-1},
\end{equation}
is alternating with terms that monotonically decrease in absolute value. Therefore, as discussed in \S \ref{part: descr, basic prop}, it can be expressed as a series of positive terms,
\begin{equation}
\arctan x = \sum_{j=1}^\infty \left( \frac{x^{4j-3}}{4j-3} - \frac{x^{4j-1}}{4j-1} \right),
\end{equation}
and the sum of the first $N$ terms has an error bounded by $x^{4N+1}/(4N+1)$. Using this into \eqref{eq: Euler Machin} yields the representation $\pi/4 = \sum_{j=1}^\infty a_j$ with
\begin{equation}
\label{eq: series, pidiv4}
a_j = \frac{2^{-4j+3} + 3^{-4j+3}}{4j-3} - \frac{2^{-4j+1} + 3^{-4j+1}}{4j-1},
\end{equation}
where $a_j$ is positive and rational; and $\pi/4 - \sum_{j=1}^N a_j < \error(N)$ for $N \geq 1$ with
\begin{equation}
\label{eq: series, pidiv4, error}
\error(N) = \frac{2^{-4N-1} + 3^{-4N-1}}{4N+1},
\end{equation}
which is rational and monotonically decreasing with $\lim_{N \rightarrow \infty} \error(N) = 0$.

\begin{figure}
\centering
\includegraphics[width = .75\textwidth]{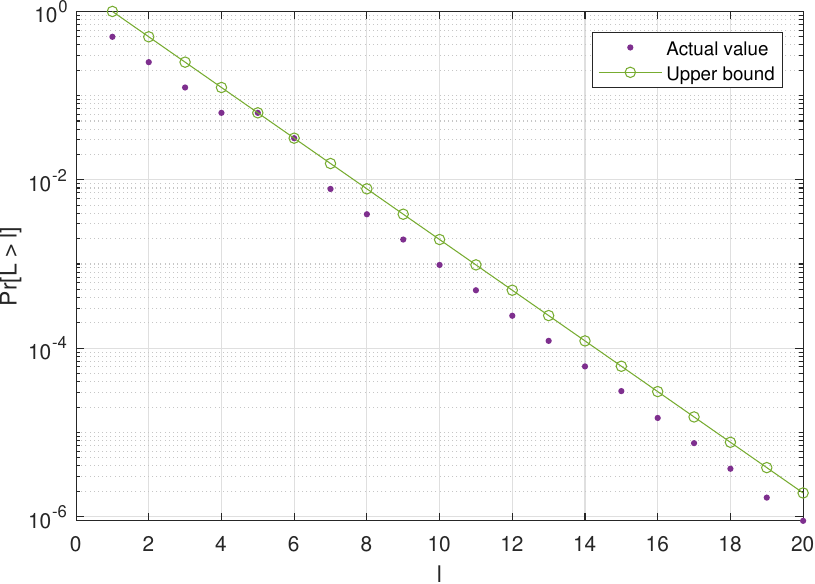}
\caption{
Sample estimate of $\Pr[L > l]$ for simulation of $\pi/4$}
\label{fig: pidiv4_PrL}
\end{figure}

\begin{figure}
\centering
\includegraphics[width = .75\textwidth]{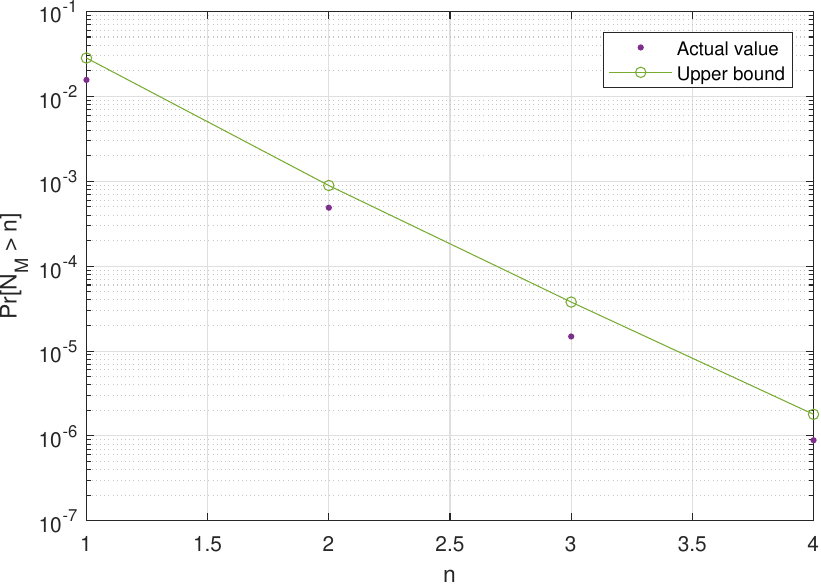}
\caption{
Sample estimate of $\Pr[N_M > n]$ for simulation of $\pi/4$}
\label{fig: pidiv4_PrNM}
\end{figure}

Based on the above, Algorithm \ref{algo: main} can be applied to \eqref{eq: series, pidiv4} and \eqref{eq: series, pidiv4, error}. The results for $10^8$ simulations are presented in Table \ref{tab: pidiv4} and in Figures \ref{fig: pidiv4_PrL} and \ref{fig: pidiv4_PrNM}. The average number of inputs is again slightly larger than the optimum $2$.

The convergence of the series is much faster in this case than in \S \ref{part: sim gamma} (exponential instead or inverse power law), which translates into smaller values of $N_M$. In fact, the maximum value observed in the $10^8$ simulations is $N_M=6$, which only occurs in $3$ cases (thus the sample estimate of $\Pr[N_M>5]$ is not reliable, and is not plotted in Figure \ref{fig: pidiv4_PrNM}). In view of \eqref{eq: series, pidiv4} and \eqref{eq: series, pidiv4, error}, the number of arithmetical operations needed for each new term of the series and for the corresponding error bound is small, and a simulation requires an average number of operations of the order of a few tens.

\begin{table}
\caption{Sample results for simulation of $\pi/4$}
\label{tab: pidiv4}
\begin{center}
  \begin{tabular}{|r|r|r|r|}
    \hline
    \multicolumn{1}{|p{2.7cm}}{\centering $\target$} & \multicolumn{1}{|p{2.7cm}}{\centering Average of $Y$} &  \multicolumn{1}{|p{2.7cm}}{\centering Average of $L$} & \multicolumn{1}{|p{2.7cm}|}{\centering Average of $N_M$} \\ \hline\hline
    $0.785398\ldots$ & $0.785402$ & $2.0467$ & $1.0161$ \\ \hline
  \end{tabular}
\end{center}
\end{table}

\subsection{Other examples}
\label{part: other examples}

The following briefly discusses how the proposed algorithm can be applied to a few other specific values of $\target$.
\begin{itemize}
\item
$\target = 1/\sqrt{2}$: using the Taylor expansion of $1/\sqrt{x}$ about $x=1$, the value $1/\sqrt{2}$ can be expressed as in \eqref{eq: serie, alt} with
\begin{equation}
b_j = \begin{cases}
1 & \text{if } j=1 \\
\displaystyle \frac{(2j-3)!!}{(2j-2)!!} & \text{if } j \geq 2.
\end{cases}
\end{equation}
Making use of the technique described in \S \ref{part: descr, basic prop}, this can be rewritten as a series with positive, rational terms, and a rational bound for the truncation error is easily obtained. Therefore Algorithm \ref{algo: main} is applicable.

\item
$\target = 1/e$: similarly, the Taylor expansion of $e^{-x}$ about $x=0$ yields an alternating series for $1/e$ of the form \eqref{eq: serie, alt} with
\begin{equation}
b_j = \frac{1}{(j-1)!},
\end{equation}
to which an analogous procedure can be applied.

\item
$\target = 1/\pi$: in this case the series for $2/\pi$ given by \cite[equation (10.1)]{Baruah09} can be employed. Adapting it for $\target = 1/\pi$, this corresponds to \eqref{eq: serie, alt} with
\begin{equation}
b_j = \begin{cases}
1/2 & \text{if } j=1 \\
\displaystyle
\frac{4j-3}{2} \left( \frac{(2j-3)!!}{(2j-2)!!} \right)^3  & \text{if } j \geq 2,
\end{cases}
\end{equation}
which is used similarly.

\item
$\target = 1/(\sqrt{2}\pi)$: a well-known series found by Ramanujan \cite[equation (2.1)]{Baruah09} can be put in the form \eqref{eq: serie, pos} with $\target = 1/(\sqrt{2}\pi)$ and
\begin{equation}
\label{eq: from Ramanujan}
a_j = 19602\, \frac{(4j-4)!}{4^{4j-4} (j-1)!^4} \, \frac{26390j-25287}{99^{4j}},
\end{equation}
where $a_j$ is rational. The first fraction in \eqref{eq: from Ramanujan} is less than $1$, and thus the truncation error can be bounded as
\begin{equation}
\label{eq: from Ramanujan, erro bound}
\sum_{j=N+1}^\infty a_j < \error(N) = 19602 \sum_{j=N+1}^\infty \frac{26390j-25287}{99^{4j}}.
\end{equation}
An explicit, rational expression for $\error(N)$ is obtained from \eqref{eq: from Ramanujan, erro bound} making use of the identities 
\begin{align}
\sum_{j=1}^\infty q^{j-1} &= \frac{1}{1-q} \\
\sum_{j=1}^\infty j q^{j-1} &= \frac{1}{(1-q)^2},
\end{align}
and then Algorithm \ref{algo: main} can be applied. 
\end{itemize}

\section{Conclusions and future work}
\label{part: concl}

An algorithm has been proposed that generates a Bernoulli random variable of arbitrary parameter $\target$, using a sequence of independent Bernoulli variables of parameter $1/2$ as input. The algorithm requires a positive series representation of $\target$, and a bound for the truncation error that converges to $0$. If the series terms and the error bound are rational, the algorithm can be implemented using only rational arithmetic.

Standard simulation methods, as implemented in computer software, typically define the parameter $\target$ as a floating-point value. This inevitably incurs a loss of precision. More specifically, it is not possible to represent irrational values (or even certain rational values) exactly using floating-point variables. The method presented in this paper avoids that problem, and generates a random variable with the exact parameter $\target$.

The algorithm consumes a random number of inputs $L$ with $2 \leq E[L] \leq 3$. Thus $\E[L]$ is close to the optimum achievable by any algorithm, which is $2$. In addition, the algorithm is fast in the sense of Nacu-Peres (that is, $\Pr[L > l]$ has an exponential bound). The number of series terms that need to be computed, $N_M$, is also random, and $\Pr[N_M > n]$ decreases with $n$ at least as fast as the truncation error bound. $\E[N_M]$ is finite if $\Pr[N_M > n]$ is $O(1/n^r)$, $r>1$.

The algorithm has been applied to the simulation of $\gamma$ and of $\pi/4$.
The former solves an open question in \cite{Flajolet11}, and establishes that $\gamma$ can be simulated using rational arithmetic only. In both cases the average numbers of consumed inputs and of required operations are very small.

As future work, an interesting extension would be to allow the inputs to be biased coins, with an arbitrary and possibly unknown parameter $p$. A straightforward approach for this problem is to transform the input $p$-coins into $1/2$-coins and then apply the algorithm presented here. With $p$ unknown, the average number of $1/2$-coins that can be obtained per $p$-coin is known to be at most $-\log_2 p - \log_2(1-p)$, and a rate arbitrarily close to this can be achieved using Peres' iterative version of von Neumann's procedure \cite{Peres92}. Consequently, the average number of $p$-coin inputs needed to generate a $\target$-coin using this approach can be roughly approximated by $\E[L] / (-\log_2 p - \log_2(1-p))$, with $\E[L]$ as resulting from Algorithm \ref{algo: main} (and bounded by Theorem \ref{theo: number of inputs}). Perhaps a more efficient method can be found.

\section{Proofs}
\label{part: proofs}

\subsection{Proof of Theorem \ref{theo: basic}}

1. Since the number of iterations $M$ is a shifted geometric random variable, it is finite with probability $1$.

2. The algorithm outputs $Y=1$ if $U < \target$ and $Y=0$ if $U > \target$, where $U$ is uniform on the interval $(0,1]$. Therefore $\Pr[Y=1] = \Pr[U < \target] = \target$ and $\Pr[Y=0] = \Pr[U > \target] = 1-\target$ (the event $U=\target$ has probability $0$).

3. Under the stated assumptions, all variables used by the algorithm are obtained as a finite sequence of additions, multiplications or divisions of rational numbers, or powers of rational numbers with integer exponents.
\qed

\subsection{Proof of Theorem \ref{theo: number of inputs}}

The number of inputs, $L$, satisfies
\begin{equation}
\label{eq: L M}
M \leq L \leq M+1.
\end{equation}
For $l \geq 1$, using \eqref{eq: Pr M geq} and \eqref{eq: L M},
\begin{equation}
\Pr[L > l] \leq \Pr[M > l-1] = \Pr[M \geq l] = 2^{-l+1}.
\end{equation}
As a consequence, \eqref{eq: fast} holds with $C=2$, $\rho= 1/2$.

Since $M$ is a shifted geometric variable with parameter $1/2$, $\E[M]=2$. Inequality \eqref{eq: bound on E L} then follows from \eqref{eq: L M}.
\qed

\subsection{Proof of Theorem \ref{theo: number of series terms}}

The sequence formed by the numbers of series terms used in each iteration, $N_k$, is deterministic, and is dictated by the values of the terms $a_j$ and of the error function $\error(n)$. This sequence is monotonically non-decreasing. The total number of terms used by the algorithm is that corresponding to the last iteration, that is, $N_M$. This is a random variable, because $M$ is.

In general, the sequence $N_k$ consists of runs of equal values, where each run has length $1$ or greater. Let $u(i)$ denote the initial index of the $i$-th run. Note that $N_{u(1)} = N_1$. Thus $N_{u(1)}, N_{u(2)}, \ldots$ is the subsequence of unique values of sequence $N_k$, with
\begin{equation}
\label{eq: subseq: incr}
N_{u(i-1)} = N_{u(i)-1} < N_{u(i)}. 
\end{equation}

The following observation is key to the proof. If the interval $(\lowerun_k, \upperun_k]$ were chosen such that $\upperun_k-\lowerun_k \leq 2^{-(k+1)}$, this would guarantee that at least one of the conditions \eqref{eq: rules for N k and s k} would hold (see Figure \ref{fig: from iteration k-1 to k}, and observe that $2^{-(k+1)}$ is $1/4$ of the length of the previous quantized interval $(\lowerqu_{k-1}, \upperqu_{k-1}]$).

Let the sequence $N_k$ and the error function $\error$ be extended by defining $N_0=0$, $\error(0)=1$. Consider $k = u(i)$ for $i \in \mathbb N$ arbitrary; that is, the index $k$ starts a run of equal values of the sequence $N_k$. This means that, at iteration $k$, the previous $N_{k-1}$ was not enough to fulfil any of conditions \eqref{eq: rules for N k and s k}, and new series terms up to index $N_k$ had to be added. In particular, $N_k-1$ terms were not sufficient to fulfil \eqref{eq: rules for N k and s k}. Taking into account the observation in the preceding paragraph, it follows that $\error(N_{k}-1)$ is necessarily greater than $2^{-(k+1)}$; that is,
\begin{equation}
\label{eq: error ineq}
\error(N_{u(i)}-1) > 2^{-(u(i)+1)}.
\end{equation}
This holds even if $i=1$ and $N_1=1$, thanks to the definitions of $N_0$ and $\error(0)$. 

Since $u(i)$ is the starting index of its run, $\Pr[N_M \geq N_{u(i)}]$ can be bounded from \eqref{eq: Pr M geq} and \eqref{eq: error ineq} as
\begin{equation}
\label{eq: Pr N M N u i}
\Pr[N_M \geq N_{u(i)}] = \Pr[M \geq u(i)] = 2^{-u(i)+1} < 4 \error(N_{u(i)}-1).
\end{equation}
The variable $N_M$ can only take the values $N_{u(1)}, N_{u(2)}, \ldots$, which form an increasing sequence. Therefore, $N_M$ cannot take any value between $N_{u(i-1)}$ and $N_{u(i)}$. Thus for $n \in \{N_{u(i-1)}, N_{u(i-1)}+1, \ldots, N_{u(i)}-1\}$ the event $N_M > n$ is equivalent to $N_M \geq N_{u(i)}$. Using the fact that $\error$ is non-increasing, \eqref{eq: Pr N M N u i} implies that, for the referred values of $n$,
\begin{equation}
\label{eq: Pr N M n}
\Pr[N_M > n] < 4 \error(N_{u(i)}-1) \leq 4 \error(n).
\end{equation}
As this is valid for $i \geq 1$, it follows that \eqref{eq: Pr N M n} holds for any $n \geq 1$.

If $\error(n)$ is $O(1/n^r)$ there exist $K$ and $n_0$ such that $\error(n) \leq K / n^r$ for all $n \geq n_0$. Therefore the mean of $N_M$ can be expressed as
\begin{equation}
\label{eq: E M N}
\begin{split}
\E[N_M] &=
\sum_{n=1}^\infty \Pr[N_M \geq n] =
\sum_{n=0}^{n_0-1} \Pr[N_M > n] + \sum_{n=n_0}^\infty \Pr[N_M > n] \\
&< \sum_{n=0}^{n_0-1} \Pr[N_M > n] + 4 \sum_{n=n_0}^\infty \error(n) \\
&\leq \sum_{n=0}^{n_0-1} \Pr[N_M > n] + 4K \sum_{n=n_0}^\infty \frac{1}{n^r}
\end{split}
\end{equation}
For $r>1$ the last term in \eqref{eq: E M N} is a convergent series \cite[corollary 2.4.7]{Abbott01}, and therefore $\E[N_M]$ is finite.
\qed

\subsection{Proof of Proposition \ref{prop: error bound for gamma}}

Noting that $\nbd(t) = \lfloor \log_2(2t) \rfloor$, the truncation error can be bounded as
\begin{equation}
\label{eq: gamma, error, 1}
\begin{split}
\gamma - \sum_{j=1}^N a_j
&= \sum_{j=N+1}^\infty \frac{\nbd(j-1)}{2j(2j-1)(2j-2)}
= \sum_{j=N}^\infty \frac{\nbd(j)}{2j(2j+1)(2j+2)} \\
&< \sum_{j=N}^\infty \frac{\log_2(2j)}{8j^3}.
\end{split}
\end{equation}
Each term $\log_2(2j) / 8j^3$ in \eqref{eq: gamma, error, 1} satisfies
\begin{equation}
\frac{\log_2(2j)}{8j^3} = \int_{j-1}^j \frac{\log_2(2j)}{8j^3} \mathrm dx < \int_{j-1}^j \frac{\log_2(2x+1)}{8x^3} \mathrm dx,
\end{equation}
and therefore
\begin{equation}\label{eq: gamma, error, 2}
\begin{split}
\gamma - \sum_{j=1}^N a_j &< \int_{N-1}^\infty \frac{\log_2(2x+1)}{8x^3} \mathrm dx < \int_{N-1}^\infty \frac{\log_2 \left( 2x \frac{2(N-1)+1}{2(N-1)} \right)}{8x^3} \mathrm dx \\
&= \int_{N-1}^\infty \frac{\log_2 x + 1 + \log_2 \left( 1+\frac{1}{2(N-1)} \right)}{8x^3} \mathrm dx.
\end{split}
\end{equation}
Using the fact that
\begin{equation}
\log_2 \left( 1+\frac{1}{2(N-1)} \right) < \frac{1}{2(N-1)\logn 2},
\end{equation}
inequality \eqref{eq: gamma, error, 2} becomes
\begin{equation}
\gamma - \sum_{j=1}^N a_j < \frac{1}{8 \logn 2} \int_{N-1}^\infty \frac{\logn x }{x^3} \mathrm dx + \left(\frac 1 8 + \frac{1}{16(N-1) \logn 2} \right) \int_{N-1}^\infty \frac{\mathrm dx}{x^3}.
\end{equation}
Evaluating the two integrals gives
\begin{equation}
\begin{split}
\gamma - \sum_{j=1}^N a_j &< \frac 1 {32(N-1)^2} \left( \frac{2\logn (N-1) + 1 + 1/(N-1)} {\logn 2} + 2 \right) \\
&= \frac 1 {16(N-1)^2} \left( \log_2 (N-1) + 1 + \frac{1 + 1/(N-1)} {2\logn 2}\right).
\end{split}
\end{equation}
Since $2\logn 2 > 1$,
\begin{equation}
\gamma - \sum_{j=1}^N a_j <
\frac {\lfloor\log_2 (N-1)\rfloor + 3 + 1/(N-1))} {16(N-1)^2} = \frac{\nbd(N-1) + 2 + 1/(N-1)}{16(N-1)^2}.
\end{equation}
\qed

\section*{Acknowledgment}

Thanks to Peter Occil for bringing to the author's attention the problem of simulating Euler's constant without using floating-point arithmetic, which led to the algorithm and results presented in this paper; also for pointing out reference \cite{Kozen14}, and for some  corrections to an early version of the manuscript.


\end{document}